\documentclass[12pt]{article}
\usepackage{amsmath}
\begin{document}
\title{An analytical exercise\footnote{Delivered October 3, 1776. Originally published as
\emph{Exercitatio analytica}, Nova Acta Academiae Scientarum Imperialis Petropolitinae \textbf{8} (1794),
69-72, and republished in \emph{Leonhard Euler, Opera Omnia}, Series 1:
Opera mathematica,
Volume 16, Birkh\"auser, 1992. A copy of the original text is available
electronically at the Euler Archive, at www.eulerarchive.org. This paper
is E664 in the Enestr\"om index.}}
\author{Leonhard Euler\footnote{Date of translation: December 8, 2004.
Translated from the Latin
by Jordan Bell, 2nd year undergraduate in Honours Mathematics, School of Mathematics and Statistics, Carleton University,
Ottawa, Ontario, Canada.
Email: jbell3@connect.carleton.ca.
Part of this translation was written
during an NSERC USRA supervised by Dr. B. Stevens.
}}
\date{}

\maketitle

1. In considering the infinite product for expressing cosine in terms
of any angle, which is
\[\cos{{\frac{\pi}{2n}}} = (1-\frac{1}{nn})(1-\frac{1}{9nn})(1-\frac{1}{25nn})
(1-\frac{1}{49nn}) \textrm{ etc.,} \]
it comes to mind to investigate a method which itself works by the nature of the
value of this product, which we in fact already know to be equal to $\cos{{\frac{\pi}{2n}}}$.
Indeed, it is possible, though with many labors which were brought about by
artifice, for a geometric demonstration of this to be given which I do not consider
at all unpleasant.

2. Therefore I set,
\[ S=(1-\frac{1}{nn})(1-\frac{1}{9nn})(1-\frac{1}{25nn}) \textrm{ etc.}, \]
and then the logarithms produced by me are summed: 
\[ lS=l(1-\frac{1}{nn})+l(1-\frac{1}{9nn})+l(1-\frac{1}{25nn})+\textrm{ etc.}, \]
and since
\[ l(1-\frac{1}{x})=-\frac{1}{x}-\frac{1}{2xx}-\frac{1}{3x^3}-\frac{1}{4x^4}
-\textrm{etc.}, \]
the successive series set out with their signs changed are:
\begin{align*}
-lS=&\frac{1}{nn}+\frac{1}{2n^4}+\frac{1}{3n^6}+\frac{1}{4n^8}+\textrm{etc.}\\
&+\frac{1}{9nn}+\frac{1}{2\cdot 9^2n^4}+\frac{1}{3\cdot 9^3n^6}+\frac{1}{4\cdot
9^4n^8}+\textrm{etc.}\\
&+\frac{1}{25nn}+\frac{1}{2\cdot 25^2n^4}+\frac{1}{3\cdot 25^3n^6}+
\frac{1}{4\cdot 25^4n^8}+\textrm{etc.}\\
&+\frac{1}{49nn}+\frac{1}{2\cdot 49^2n^4}+\frac{1}{3\cdot 49^3n^6}+\frac{1}{4\cdot
49^4n^8}+\textrm{etc.}\\
&\textrm{etc.}
\end{align*}

3. But now if we arrange each of the rows vertically in columns, we
will obtain the following series for $-lS$:
\begin{align*}
-lS=&\frac{1}{nn}(1+\frac{1}{3^2}+\frac{1}{5^2}+\frac{1}{7^2}+\frac{1}{9^2}+
\textrm{etc.})\\
&+\frac{1}{2n^4}(1+\frac{1}{3^4}+\frac{1}{5^4}+\frac{1}{7^4}+\frac{1}{9^4}+\textrm{etc.})\\
&+\frac{1}{3n^6}(1+\frac{1}{3^6}+\frac{1}{5^6}+\frac{1}{7^6}+\frac{1}{9^6}+
\textrm{etc.}\\
&+\frac{1}{4n^8}(1+\frac{1}{3^8}+\frac{1}{5^8}+\frac{1}{7^8}+\frac{1}{9^8}+
\textrm{etc.}\\
&\textrm{etc.}
\end{align*}
Thus in this way, the investigation is led to the summation of series of
even powers of the harmonic progression, $1, \frac{1}{3},\frac{1}{5},
\frac{1}{7}$, etc.

4. By putting $\frac{\pi}{2}=\rho$ for the sake of brevity, it can
then be seen that the sums of powers can be represented in the following
manner:
\begin{align*}
1+\frac{1}{3^2}+\frac{1}{5^2}+\frac{1}{7^2}+\textrm{etc.}=A\rho^2,\\
1+\frac{1}{3^4}+\frac{1}{5^4}+\frac{1}{7^4}+\textrm{etc.}=B\rho^4,\\
1+\frac{1}{3^6}+\frac{1}{5^6}+\frac{1}{7^6}+\textrm{etc.}=C\rho^6,\\
\textrm{etc.}
\end{align*}
with the first one taken as $A=\frac{1}{2}$, and then the other letters determined
by the preceding ones in the following way:
\begin{align*}
B=\frac{2}{3}A^2, C=\frac{2}{5}\cdot 2AB, D=\frac{2}{7}(2AC+BB),\\
E=\frac{2}{9}(2AD+2BC), F=\frac{2}{11}(2AE+2BD+CC),\textrm{etc.}
\end{align*}
the truth of which will soon shine forth from a most beautiful conjunction
of analysis.

5. With the substituted values, this series is obtained:
\[ -lS=\frac{a\rho^2}{nn}+\frac{1}{2}\cdot \frac{B\rho^4}{n^4}+
\frac{1}{3}\cdot \frac{C\rho^6}{n^6}+\frac{1}{4}\cdot \frac{D\rho^8}{n^8}+
\textrm{etc.} \]
If we put $\frac{\rho}{n}=x$, so that $x=\frac{\pi}{2n}$, then this
series is obtained:
\[ -lS=Axx+\frac{1}{2}Bx^4+\frac{1}{3}Cx^6+\frac{1}{4}Dx^8+\textrm{etc.} \]
So that we can get rid of the fractions $\frac{1}{2},\frac{1}{3},\frac{1}{4}$,
etc., we differentiate, and then after dividing the result by $2\partial x$
it follows,
\[ -\frac{\partial S}{2S\partial x}=Ax+Bx^3+Cx^5+Dx^7+\textrm{etc.} \]

6. Then for the sake of brevity we set $-\frac{\partial S}{2S\partial x}=t$,
so that we have
\[ t=Ax+Bx^3+Cx^5+Dx^7 + \textrm{etc.}, \]
for which the square is seen to be this series:
\begin{align*}
tt=&A^2xx&+2ABx^4&+2ACx^6&+2ADx^8&+2AEx^{10}&+\textrm{etc.}\\
&&&+BB&+2BC&+2BD&+\textrm{etc.}\\
&&&&&+CC&+\textrm{etc.}
\end{align*}
and in this, each of the powers of $x$ are the very same as those found
in the formulas which express the determination of the letters $A,B,C,D$:
the only thing missing are the coefficients $\frac{2}{3},\frac{2}{5},
\frac{2}{7}$, etc.

7. But we can indeed introduce these coefficients by integrating after
we have multiplied through by $2\partial x$. So it will be obtained that
\begin{align*}
2\int tt\partial x=\frac{2}{3}A^2x^3+\frac{2}{3}\cdot 2ABx^5+\frac{2}{7}
(2AC+BB)x^7+\\
\frac{2}{9}(2AD+2BC)x^9+\frac{9}{11}(2AE+2BD+CC)x^{11}+
\textrm{etc.}
\end{align*}
Since at present,
\[ \frac{2}{3}A^2=B,\frac{2}{3}\cdot 2AB=C,\frac{2}{7}(2AC+BB)=D,\textrm{etc.}, \]
with these values substituted in, we then attain this series:
\[ 2\int tt\partial x=Bx^3+Cx^5+Dx^7+Ex^9+\textrm{etc.} \]

8. Therefore with the series we had from before:
\[ t=Ax+Bx^3+Cx^5+Dx^7+\textrm{etc.}, \]
this equation clearly follows:
\[ t=Ax+2\int tt\partial x, \]
which when differentiated yields
\[ \partial t=A\partial x+2tt\partial x=\frac{1}{2}\partial x+2tt\partial x,
\textrm{ because } A=\frac{1}{2}. \]
From this we therefore have $2\partial t=\partial x(1+4tt)$, from which it
is $\partial x=\frac{2\partial t}{1+4tt}$, whose integral is
at hand, namely $x=\textrm{A tang.} 2t$,\footnote{Translator:
Euler uses here A tang. $2t$ to express the arctangent of $2t$, and later
uses tang. $x$ to stand for the tangent of $x$.}
for which the addition of a constant is not
desired, since by setting $x=0$ then $t$ should simultaneously vanish. 
Hence with this equation that has been discovered, if the quantity
$x$ is seen as an angle then it will be $2t=\textrm{tang.} x$.
Since it was
already $t=-\frac{\partial S}{2S\partial x}$, from this are deduced
these equations:
\[ -\frac{\partial S}{S\partial x}=\textrm{tang.} x,
 \textrm{ and } -\frac{\partial S}{S}=
\frac{\partial x \sin{x}}{\cos{x}}. \]

9. Therefore, since it is $\partial x \sin{x}=-\partial \cos{x}$, it will
be $\frac{\partial S}{S}=\frac{\partial \cos{x}}{\cos{x}}$, and then
by integrating, $lS=l\cos{x} +C$, which to be consistent should be defined such
that setting $x=0$ makes $lS=0$. Thus it will be $C=0$, so
$lS=l\cos{x}$, and then so that the numbers can come forth, $S=\cos{x}$.

10. Moreover, since we already had that $x=\frac{\pi}{2n}$, from this
it is clear that the value being sought out for $S$ is revealed to be
$S=\cos{\frac{\pi}{2n}}$, which agrees totally with what was known from before.
Therefore, analysis most excellently confirms the former relation between the
letters $A,B,C,D$, which I introduced in earlier calculations.

\end{document}